\documentclass[12pt]{article}
\usepackage{amsmath, amssymb, eucal, latexsym}
\usepackage[cp1251]{inputenc}
\usepackage[english]{babel}
\usepackage{graphicx}

\def\N{\Bbb N}
\def\R{\Bbb R}
\def\F{\Bbb F}
\def\P{\Bbb P}

\def\H{\mathsf H}

\def\E{\mathsf {E}}
\def\t{\tilde}

\def\eps{\varepsilon}

\def\Span{{\rm Span\,}}
\def\Cay{{\mathbf{ Cay}}}
\def\Gr{{\mathbf G\,}}

\newtheorem{theorem}{Theorem}
\newtheorem{proposition}[theorem]{Proposition}
\newtheorem{lemma}[theorem]{Lemma}
\newtheorem{cor}[theorem]{Corollary}

\begin{document}

 \begin{center}
\textbf{ON SUBGRAPHS OF RANDOM CAYLEY SUM GRAPHS}
\end{center}

 \begin{center}
                                                         S. V. KONYAGIN, I. D. SHKREDOV
\footnote{
This work is supported by the Russian Science Foundation under grant 14-50-00005.
}\\

    \end{center}

\bigskip

\begin{center}
    Abstract.
\end{center}

{\it \small
We prove that asymptotically almost surely, the random Cayley
sum graph over a finite
%IS
%elementary abelian 2--group
abelian group
$\Gr$ has edge
density close to the expected one on every induced subgraph of
size at least $\log^c |\Gr|$, for any fixed $c > 1$ and $|\Gr|$ large enough.
}

\bigskip
\section{Introduction}
\label{sec:introduction}
\bigskip

Let $A$ be a subset of an additively written group $\Gr$.
We denote by $\Cay (A,\mathbf{G})$ the {\it Cayley sum graph} induced by $A$ on $\Gr$, which
is the directed graph on the vertex set $\Gr$ in which $(x,y) \in \Gr \times \Gr$ is an edge if and only if
$x+y\in A$ ($x=y$ is allowed).
Such graphs are classical combinatorial objects, see, e.g. \cite{H_B}.
B. Green \cite{Green} initiated to study  the {\it random
Cayley sum graph}, considering finite groups $\Gr$ and selecting
$A$ at random by choosing each $x\in \mathbf{G}$ to lie in $A$
independently and at random with probability $1/2$. General random
graphs are considered in \cite{Bollobas}.
Results about random Cayley sum graphs can be found, for example,
in \cite{Green}, \cite{GM}, \cite{M}, \cite{M2}.
R. Mrazovi\'{c} \cite{M}  proved the following theorem.

\begin{theorem}
Let $\mathbf{G}$ be a
finite group
and $w: \mathbb{N} \to
\mathbb{R}$
be a growing function that tends to infinity. Let $A \subset
\mathbf{G}$ be a random subset obtained by putting every element
of $\mathbf{G}$ into $A$ independently with probability
$\frac{1}{2}$. Then with probability $1-o(1)$, for all sets
    $X,Y\subset \mathbf{G}$
    %,
    with
    $$
        |X|\ge w(|\Gr|) \log |\Gr| \quad  \mbox{ and }  \quad |Y| \ge w(|\Gr|) \log^2 |\Gr| %\,,
    $$
    one has
    \begin{equation}\label{f:Mrazovic}
        \sum_{x\in X} \sum_{y\in Y} A(x+y) = \frac{1}{2} |X| |Y| + o(|X| |Y|) \,,
        \quad (|\Gr|\to\infty)\,,
    \end{equation}
    where the rate of convergence implied by the $o$--notation depends only on $w$.
\label{t:Mrazovic}
\end{theorem}

In our paper for a set $A$ we  use the same letter to denote its
characteristic function $A : \Gr \to  \{ 0, 1\}$.

In the same paper Mrazovi\'{c} showed that
there is no $C$ such that the assumption of Theorem \ref{t:Mrazovic} can be relaxed to
 $\min \{ |X|, |Y|\} \ge C\log |\Gr| \log \log |\Gr|$.

Theorem \ref{t:Mrazovic} shows that with high probability, the edge density of the random Cayley sum
graph on all induced subgraphs of size at least $\log^{2+\eps} |\Gr|$ is close to $1/2$.

Using some tools from Additive Combinatorics, we show that Theorem \ref{t:Mrazovic}
can be improved.
%IS
% for some abelian groups.

\begin{theorem}
	Let ${\Gr}$ be a finite abelian group of size $N$
	and $w: \mathbb{N} \to \mathbb{R}$ be
	a growing function that tends to infinity.
	Let $A \subset \mathbf{G}$ be a random subset obtained by putting every
	element of $\mathbf{G}$ into $A$ independently with probability
	$\frac{1}{2}$. Then with probability
	$1-o(1)$, for all sets $X,Y\subset \mathbf{G}$ such that
	$$|X|\ge w(|\Gr|)\log |\Gr| (\log \log |\Gr|)^2,\quad
	|Y| \ge w(|\Gr|) \log |\Gr|(\log \log |\Gr|)^{10} \,,
	$$
	one has
	$$
	\sum_{x\in X} \sum_{y\in Y} A(x+y) = \frac{1}{2} |X| |Y| + o(|X| |Y|)
	\quad (|\Gr|\to\infty)\,,
	$$
	where the rate of convergence implied by the $o$--notation depends only on $w$.
	\label{t:main3}
\end{theorem}

Thus lower and upper bounds for size of sets $X,Y$ differ by some powers of double logarithms. 

Let us say a few words about the proof.

It was showed in \cite{M} that if for some $X,Y$ the sum of the left-hand
side of (\ref{f:Mrazovic}) deviates significantly from $\frac12|X||Y|$,
then  the common energy (see the definition in the next section) of $X$ and $Y$ must be close
to the trivial upper bound $|X||Y| \min \{ |X|, |Y|\}$.
Mrazovi\'{c} used a random choice to avoid such
%IS (и ниже)
a
situation (see details in \cite{M}) and using structural results from \cite{SS_dim}, \cite{SY} we add one more twist to his arguments,
hence proving that large portions of $X,Y$ must be very structured in this case.
%%%New I.S.
It follows that the number of such sets is much smaller
than the number of all possible pairs of arbitrary sets $X,Y$.
This allows us to relax  the conditions on sizes of $|X|,|Y|$ and to obtain our bound
$\log^{c} |\Gr|$ with $c>1$.

%IS
First we consider the case of elementary abelian 2--groups and prove Theorem \ref{t:main3}
%SK in this situation.
in this situation (with $c>3/2$) using some arguments from \cite{M}.
For such groups the proof is simpler and more transparent.
For general case see sections 
%\ref{sec:estim}, \ref{sec:main}.
\ref{sec:estim}, \ref{sec:improvement}. 

\bigskip

We thank Rudi  Mrazovi\'{c}, Mikhail Gabdullin for fruitful discussions, 
%for careful 
and the reviewers for their useful remarks.

\bigskip
\section{Definitions and preliminary results}
\label{sec:definitions}
\bigskip

%SK
%Although Theorem~\ref{t:main} concerns with finite groups of exponent $2$ only,
%in this section we allow $\mathbf{G}$ to be any finite abelian group.

Let $\Gr$ be an abelian group. The {\it additive energy $\E
(A,B)$} between two sets $A$ and $B$ from $\Gr$ is
(see \cite{TV})
$$
    \E (A,B) = |\{ (a_1,a_2,b_1,b_2) \in A\times A \times B\times B ~:~ a_1+b_1 = a_2+b_2  \}| \,.
$$
The {\it sumset} of $A$ and $B$ is
$$
    A+B := \{ a+b ~:~ a\in A,\, b\in B \} \,.
$$
By $A\bigsqcup B$ denote the  union of two disjoint sets $A,B$.

\bigskip

Recall  a simple lemma, see, e.g., \cite[Lemma 12]{S_exp_p}.

\begin{lemma}
    For any finite sets $X,Y,Z\subset \Gr$ one has
$$
    \E (X\cup Y,Z)^{1/2} \le \E (X,Z)^{1/2} + \E (Y,Z)^{1/2} \,,
$$
    and for disjoint union of $X$ and $Y$ the following holds
$$
    \E (X\sqcup Y,Z) \ge \E (X,Z) + \E (Y,Z) \,.
$$
\label{l:norm_AB}
\end{lemma}

Now let us recall the notion of the (additive) {\it dimension} of
a set. A finite set $\Lambda \subset \mathbf{G}$ is called {\it
dissociated} if any equality of the form
$$\sum_{\lambda \in \Lambda}\eps_{\lambda}\lambda=0$$
for $\eps_{\lambda}\in \{-1,0,1\}$ implies $\eps_{\lambda}=0$ for
all $\lambda \in \Lambda.$ The notion of dissociativity appears
naturally in analysis, see \cite{Rudin_book}. The size of a
largest dissociated subset of $A$  is  called the (additive) {\it
dimension} of the set $A$ and is denoted by $\dim (A)$. For a
subset $S = \{ s_1, \dots, s_l\} \subset \mathbf{G}$ one can
define
$$
    \Span (S) := \Big\{ \sum_{j=1}^{l} \eps_j s_j ~:~ \eps_j \in \{ 0,-1,1 \} \Big\} \,.
$$
It is easily seen that if $S$ is a dissociated subset of $A$ of
size $|S| = \dim A$, then $A\subset \Span (S)$.

Notice that if $\mathbf{G}$ is a finite group of exponent $2$ (hence, a linear space over
$\F_2$), then $\Span(S)$ is the linear span of $S$, and $\dim S$ is its dimension.

\bigskip

We need the main result from \cite{SY}, also see \cite[Theorem 19]{SS_dim}.

\begin{theorem}
    Let  $A,B$ be finite non-empty subsets of
    %finite
    an
    abelian group, $|A| \ge |B|$.
    If $\E (A,B) \ge \frac{|A| |B|^2}{K}$, then there exist a non-empty set
    $B_* \subset B$ such that
    \begin{equation}\label{f:dim_sy_new}
        \dim (B_*) \ll K \log |A|  \,,
    \end{equation}
    and
    \begin{equation}\label{f:E(A,B_1)'}
        \E (A,B_*) \ge 2^{-5} \E(A,B) \,.
    \end{equation}
\label{t:SS_dim}
\end{theorem}

Theorem \ref{t:SS_dim} shows that if $\E (A,B)$ is large, then $B$
contains a large, well-structured subset $B_*$.

Let us derive a simple
consequence of the theorem above.

\begin{cor}
    Let  $A,B$ be finite subsets of
    an
    abelian group with $|A| \ge |B| \ge 2$.
    Suppose that $\E (A,B) = \frac{|A| |B|^2}{K}$, and $M\ge K$ be a parameter.
    Then
there is a partition $B=B' \bigsqcup B''$
    such that
    \begin{equation}\label{f:SS_dim1}
        \dim (B') \ll
            M \log |A| \cdot
            \log (|B|M/K)
            \,, \quad \quad \E(A,B') \gg \E(A,B) \,,
    \end{equation}
    and
    \begin{equation}\label{f:SS_dim2}
        \E(A,B'') \le \frac{|A| |B''|^2}{M} \,.
    \end{equation}
\label{c:SS_dim}
\end{cor}

{\bf Proof.}
    Our arguments is a sort of an algorithm similar to that found in \cite{KS, Sh}.
    We construct an increasing sequence of sets
    $\emptyset = B'_1  \subset B'_2 \subset \dots \subset B'_k \subset B$
    and a decreasing sequence of sets $B=B''_1 \supset B''_2 \supset \dots \supset B''_k$
    such that
    for any $j=1,2,\dots, k$ the sets $B'_j$ and $B''_{j}$  are disjoint and moreover $B = B'_j \sqcup B''_{j}$.
%SKG      If at some step $j$ we have either $\E(A,B''_j) < |A| |B''_j|^2 / M$ or $B''_j = \emptyset$
     If at some step $j$ we have either $\E(A,B''_j) < |A| |B''_j|^2 / M$ or $B''_j = \emptyset$
(notice that due to the definition of $K$ and the supposition $M\ge K$ this can happen only
for $j>1$)
then we stop our algorithm putting
    $B'' = B''_j$, $B' = B'_{j}$, and $k=j$.
    In the opposite situation where $\E(A,B''_j) \ge |A| |B''_j|^2 / M$ we apply Theorem \ref{t:SS_dim} to the set $B''_j$,
     finding a non-empty subset $G_j$ of $B''_j$ such that
     \begin{equation}\label{tmp:07.10.2016_1}
        \dim (G_j) \ll M  \log |A|  \,,
     \end{equation}
     and
     \begin{equation}\label{tmp:07.10.2016_2}
        \E (A,G_j) \ge 2^{-5} \E(A,B''_j) \,.
     \end{equation}
     After that we put $B''_{j+1} = B''_j \setminus G_j$, $B'_{j+1} = B'_{j} \sqcup G_j$ and repeat the procedure.
     Clearly, $B'_{k} = \bigsqcup_{j=1}^{k} G_j$.
In view of Lemma \ref{l:norm_AB} and (\ref{tmp:07.10.2016_2}), we get
     $$
        \E (A,B''_j) \ge \E (A,G_j) + \E (A,B''_{j+1}) \ge 2^{-5} \E(A,B''_j) + \E (A,B''_{j+1})
     $$
     whence $\E (A,B''_{j+1}) \le \frac{31}{32} \E (A,B''_j)$.
     It follows that our algorithm stops after at most
     $k\ll \log (|B|M/K)$
     steps.
    Because $G_1\subset B'_{j}$, $j\ge 2$, we have in view of (\ref{tmp:07.10.2016_2}) that for $j\ge 2$ one has
    $$
        \E(A,B'_{j}) \ge \E(A,G_1) \ge 2^{-5} \E(A,B''_1) = 2^{-5} \E(A,B)
    $$
%SKG    and thus inequality  $\E(A,B') \gg \E(A,B)$ takes place.
    and thus inequality  $\E(A,B') \gg \E(A,B)$ holds.
     Finally, from estimate (\ref{tmp:07.10.2016_1}), we obtain
$$
    \dim (B') \le \sum_{j=1}^{k-1} \dim (G_j)
        \ll kM \log |A| \ll
        M \log |A| \cdot \log (|B|M/K) \,.
$$
This completes the proof of the corollary.
$\hfill\Box$

%IS Done

\bigskip

We finish this section with a result on the number of sets with small dimension.

\begin{lemma}
    Let $\mathbf{G}$ be a finite abelian group, and write $N = |\mathbf{G}|$.
    Let $n,d\in\N$ with $n\ge2\log N$. Then the number of sets $X\subset\Gr$
    with $0<|X|\le n$ and $\dim X\le d$ is at most $e^{2nd}$.
    \label{l:count_sm_dim}
\end{lemma}

{\bf Proof.}
Take $X\subset X'=\Span(\Lambda)$ where $|\Lambda|=\dim X\le d$.
The number of sets $\Lambda$ is at most $N^d$. For a fixed $\Lambda$,
we have $|X'|\le 3^d$, and the number of sets
%IS
$X \subset X'$
(with fixed $\Lambda$)
is at most
$|X'|^n\le 3^{nd}$. Therefore, the total number of sets $X$ is
at most
$$N^d3^{nd} < e^{(\log N+1.1n)d} \le e^{2nd},$$
as required.
$\hfill\Box$

\bigskip
%IS
%\section{The proof of the main result}
\section{A model case}
\label{sec:model}
\bigskip

The main result of this section is the following

\begin{theorem}
    Let ${\Gr}$ be a finite group of exponent $2$
    and $w: \mathbb{N} \to \mathbb{R}$ be
    a growing function that tends to infinity.
    Let $A \subset \mathbf{G}$ be a random subset obtained by putting every
    element of $\mathbf{G}$ into $A$ independently with probability
    $\frac{1}{2}$. Then with probability
    $1-o(1)$, for all sets $X,Y\subset \mathbf{G}$ such that
    $$
    |X|, |Y| \ge
    w(|\Gr|) (\log \log |\Gr|) \cdot \log^{3/2} |\Gr| \,,
    $$
    one has
    $$
    \sum_{x\in X} \sum_{y\in Y} A(x+y) = \frac{1}{2} |X| |Y| + o(|X| |Y|)
    \quad (|\Gr|\to\infty)\,,
    $$
    where the rate of convergence implied by the $o$--notation depends only on $w$.
    \label{t:main}
\end{theorem}

We notice that the groups considered in Theorem \ref {t:main}
are abelian and they can be treated as vector
spaces over the field $\F_2$.

\bigskip

For finite, non--empty subsets $X,Y$, and $A$ of $\Gr$, let
$$
    \sigma_A (X,Y) := \frac{1}{|X||Y|} \sum_{x\in X} \sum_{y\in Y} A(x+y) - \frac{1}{2}
     =
        \frac{1}{|X||Y|} \sum_{x\in X} \sum_{y\in Y} \left( A(x+y) - \frac{1}{2} \right) \,.
$$

%SKG The following technical result is in the heart of \cite{M} (see section 4 of that paper).
The following technical result is the heart of \cite{M} (see section 4 of that paper).

\begin{proposition}
    Let $\mathbf{G}$ be a finite abelian group, and write $N = |\mathbf{G}|$.
    If $A$ is a random subset of $\mathbf{G}$ obtained by putting every element of
    $\Gr$ into $A$ independently with probability $\frac{1}{2}$, then for any
    $r,K\ge1$ and $\eps \in (0,1]$, the probability that there exist
    $X,Y\subset \Gr$ satisfying
 $$
    |Y| \ge |X| \ge 2000 \eps^{-4} \log N,\, |Y| \ge r,\,
    \E(X,Y) \le \frac{|X|^2 |Y|}{K},
$$
    and
$$
    |\sigma_A (X,Y)| \ge \eps
$$
is at most
\begin{equation}\label{f:p_M}
    C\exp\left( \frac{2000 \log^2 N}{\eps^4} - \frac{\eps^2 rK}{40} \right)
\end{equation}
with an absolute constant $C$.
\label{p:M}
\end{proposition}

Proposition \ref{p:M} was not stated in \cite{M} explicitly, and
for completeness, we prove it in Appendix.

\bigskip

%SKG Let us show quickly how Proposition \ref{p:M} implies Theorem \ref{t:Mrazovic}.
Let us show quickly how Proposition \ref{p:M} implies Theorem \ref{t:Mrazovic}
for abelian groups $\Gr$.
Choosing
$K=1$ and $r = C\eps^{-6} \log^2 N$ with $C$ large enough,
we obtain
%SK that the probability of existence $X,Y$, $|Y|\ge |X| \ge  2000 \eps^{-4} \log N$, $|Y| \ge r$
that the probability of existence $X,Y$, $|Y|\ge |X| \ge  2000 \eps^{-4} \log|\Gr|$, $|Y| \ge r$
such that  $|\sigma_A (X,Y)| \ge \eps$ is less than
$$
%SK    C_1 \exp(-C_2 \eps^{-4} \log^2 N ) = o(1)
%SK        \quad \mbox{ as } \quad N\to +\infty \,,
    C_1 \exp(-C_2 \eps^{-4} \log^2 |\Gr| ) = o(1)
        \quad \mbox{ as } \quad |\Gr|\to +\infty \,,
$$
where $C_1,C_2>0$ are some absolute constants.

\bigskip

The next
corollary immediately follows from Proposition \ref{p:M}
(as applied with $K=M$ and  $r=(\eps/4) w(N) (\log \log N) \cdot \log^{3/2} N$).

\begin{cor}
        Let $\mathbf{G}$ be a finite abelian group, and write $N = |\mathbf{G}|$.
    If $A$ is a random subset of $\mathbf{G}$ obtained by putting every element of
    $\Gr$ into $A$ independently with probability $\frac{1}{2}$, then for
$M = (\log\log N)^{-1} (\log N)^{1/2}$, any $\eps \in (0,1]$
such that
$$
\eps^{7} \ge \frac{2^{25}}{w(N) \log \log N \sqrt{\log N}},
$$
and any growing function $w: \mathbb{N} \to \mathbb{R}$ tending to infinity,
the probability that there exist $X,Y\subset \Gr$ satisfying
$$
|Y| \ge |X| \ge (\eps/4) w(N) (\log \log N) \cdot \log^{3/2} N,\,
\E(X,Y) \le \frac{|X|^2 |Y|}{M},
$$
and
$$
|\sigma_A (X,Y)| \ge \eps/2,
$$
tends to $0$ as $N\to\infty$.
\label {c:M}
\end{cor}

\bigskip

{\bf Proof of Theorem \ref{t:main}.}
Take a random set $A$ and suppose that for some $X,Y$ one has $|\sigma_A (X,Y)| \ge \eps$.
Without loss  of generality, suppose that $|Y| \ge |X| \ge (\eps/4) w(N) (\log \log N) \cdot \log^{3/2} N$.
In view of Theorem \ref{t:Mrazovic}, we can assume that $|X|, |Y| \ \ll \log^{5/2} N$, say.
Otherwise the probability of  the event $|\sigma_A (X,Y)| \ge \eps$ is $o(1)$.

Denote
$$K = \frac{|X|^2 |Y|}{\E(X,Y)},\quad M = (\log\log N)^{-1} (\log N)^{1/2}.$$
If $M\le K$ then we can apply Corollary~\ref{c:M} and conclude
that the probability of this event is $o(1)$. Thus we can assume
that $M\ge K$. Applying Corollary \ref{c:SS_dim} to the sets $X$,
$Y$, we find $X',X''\subset X$ such that $X=X' \bigsqcup X''$,
$\dim(X') \ll M (\log \log N)^2$, $\E(X',Y) \gg \E(X,Y)$ and
$\E(X'',Y)\le |Y| |X''|^2/M$.

We can assume that with high probability
the following holds:
\begin{equation}\label{f:supposX'}
|\sigma_A (X',Y)| \ge \eps/2,\quad |X'| \ge (\eps/2)|X|.
\end{equation}
Indeed, if one of these two inequalities does not hold, then
$$
    \left| \sum_{x\in X',\, y\in Y} \left( A(x+y) - \frac{1}{2} \right) \right| \le \eps |X||Y| /2
$$
and we have
\begin{equation*}\label{f:calc_1}
    \eps |X| |Y| \le |\sigma_A (X,Y)| |X||Y|
= \left| \sum_{x\in X,\, y\in Y} \left( A(x+y) - \frac{1}{2} \right) \right|
    \le
\end{equation*}
\begin{equation*}\label{f:calc_1.5}
    \le
     \left| \sum_{x\in X',\, y\in Y} \left( A(x+y) - \frac{1}{2} \right) \right|
    +
     \left| \sum_{x\in X'',\, y\in Y} \left( A(x+y) - \frac{1}{2} \right) \right|
        \le
\end{equation*}
\begin{equation}\label{f:calc_2}
        \le
            \eps |X||Y| /2
                +
                    |X''| |Y| |\sigma_A (X'', Y)| \,.
\end{equation}
Whence
$$
    \frac{\eps |X||Y|}{2}  \le |X''||Y| |\sigma_A (X'',Y)| \,.
$$
The last bound implies that $|X''| \ge \eps |X|/2$ and $|\sigma_A (X'',Y)| \ge \eps/2$.
Using Corollary
\ref{c:M}
with the sets $X''$, $Y$,
we see that the probability of the last inequality is $o(1)$.
Thus, we will assume
that (\ref{f:supposX'}) holds.

Split the set $Y$ onto sets $\tilde{Y}_j$ (using lexicographical  ordering, say) such that
$|X'|/2 \le |\t{Y}_j| \le |X'|$.
Then arguing
as in (\ref{f:calc_2}), we see that for some $\t{Y}_j$ one has
$\sigma_{A} (X',\t{Y}_j) \ge \eps/2$.
Indeed
$$
    2^{-1} \eps\cdot |X'| \sum_j |\t{Y}_j|
    = 2^{-1} \eps \cdot |X'| |Y| \le |\sigma_A (X',Y)| |X'| |Y|
        \le
$$
$$
        \le
        \sum_{j} \left|\sum_{x\in X',\, y\in \t{Y}_j} \left( A(x+y) - \frac{1}{2} \right) \right|
            =
                |X'| \sum_{j} \sigma_A (X',\t{Y}_j) |\t{Y}_j|
                %\,.
$$
and thus there is $j$ with $\sigma_{A} (X',\t{Y}_j) \ge \eps/2$.
Put $\t{Y} = \t{Y}_j$.
After that taking into account (\ref{f:supposX'}) and
applying  
%Proposition \ref{p:M} 
Corollary \ref{c:SS_dim} 
to the sets $X'$, $\t{Y}$, we find $Y',Y''\subset \t{Y}$
such that $\t{Y}=Y' \bigsqcup Y''$,
$\dim(Y') \ll M (\log \log N)^2$,
$\E(X',Y') \gg \E(X',\t{Y})$
and $\E(X',Y'')\le |X'| |Y''|^2/M$ (if  $\E(X',\t{Y}) := |X'||\t{Y}|^2/K < |X'| |\t{Y}|^2/M$, then it is nothing to prove, just put $Y'' = \t{Y}$, 
$Y' = \emptyset$, otherwise $M\ge K$). 
Again, we can assume with probability $1-o(1)$
\begin{equation*}\label{f:supposY'}
|\sigma_A (X',Y')| \ge \eps/4,\quad |Y'| \ge (\eps/4)|\t{Y}|.
\end{equation*}

For fixed $\eps$ and large $N$ in view of the assumption $|X|\ge w(N) (\log \log N)\cdot \log^{3/2} N$, we have
\begin{equation*}\label{f:min1}
    \min\{ \eps |\t{Y}|, |X'| \} \gg \eps^{-4} \log N \,.
\end{equation*}

Up to this point we did not use a specific structure of the group $\Gr$.
Notice that in this group for any subset $A$ the additive dimension $\dim (A)$ is just the ordinary  dimension of its  linear span. 
Consider the set  $L := \Span(Y'\cup X')$ of dimension
$$
    \dim (L) \le \dim(X')+\dim (Y')
    :=d \ll M (\log \log N)^2   \,.
$$
Recall now that $\Gr$ is a linear space over $\F_2$. 
Therefore,
$L$ is also an abelian group, and we can apply to $L$ Proposition \ref{p:M}
implying that the probability of the
inequality $\sigma_A (X',Y') \ge \eps/4$
is less than (also, see the calculations after the Proposition)
\begin{equation*}
    C' \exp\left( \frac{2^{19} d^2}{\eps^4} - \frac{\eps^2 \max\{|X'|,|Y'|\}}{160} \right) \,,
\end{equation*}
where $C'>0$ is some absolute constant.
Because the number of sets $L$ is roughly bounded by $N^{d}$,
and the number of sets $\t{Y}_j$ is bounded by $O(|Y|/|X'|)$,
we obtain that the total probability tends to zero, if
\begin{equation}\label{f:suppos_d}
    \eps^2 \max \{ |X'|, |Y'| \}
    %\gg \frac{M^2 (\log \log N)^4}{\eps^4} + M \log N (\log \log N)^2
        \gg
    \frac{d^2}{\eps^4} + d \log N + \log |Y|
        \gg
            \frac{d^2}{\eps^4} + d \log N
        \,.
\end{equation}
Since due to (\ref{f:supposX'})
$$d\ll (\log\log N) (\log N)^{1/2},\quad |X'| \ge
(\eps/2) w(N) (\log \log N) \cdot \log^{3/2} N \,,$$
we see that (\ref{f:suppos_d}) holds for large $N$.
This completes the proof.

\bigskip
%IS \section{Estimates for probabilities of large deviations}
\section{On large deviations}
\label{sec:estim}
\bigskip

%IS
%As in Section \ref{sec:model} we denote
We use the notations of section  \ref{sec:model}.
Recall, that for finite, non--empty subsets $X,Y$, and $A$ of $\Gr$,
%let
we denote
$$
    \sigma_A (X,Y) := \frac{1}{|X||Y|} \sum_{x\in X} \sum_{y\in Y} A(x+y) - \frac{1}{2}
     =
        \frac{1}{|X||Y|} \sum_{x\in X} \sum_{y\in Y} \left( A(x+y) - \frac{1}{2} \right) \,.
$$
In this section we fix $X\subset\Gr$ and estimate the probabilities that
the deviations $|\sigma_A (X,Y)|$ are large where $Y\subset\Gr$
%SK (see precise statements below).
(see precise statements below). If $Y=\{y\}$ we will write for simplicity
$\sigma_A (X,y)$ rather than $\sigma_A (X,\{y\})$.

\begin{lemma}
    Let $\mathbf{G}$ be a finite abelian group. %, and write $N = |\mathbf{G}|$.
    Fix $X\subset\Gr$ with $|X|=n$. Let $\eps \in (0,1/2]$ and $y_1,\dots,y_k\in\Gr$
    satisfy the condition
  \begin{equation}
%\label{disc2}
\label{small_inters}
\left|(X+y_i)\cap\left(\bigcup_{j=1}^{i-1}(X+y_j)\right)\right|\le\eps n
    \quad(i=2,\dots,k).
\end{equation}
%   \sum_{j=1}^{i-1}\left(|X+{y_i}\cap X+{y_j}|\right) \le \eps n
%   \quad(i=2,\dots,k).
%  \end{equation}
    If $A$ is a random subset of $\mathbf{G}$ obtained by putting every element of
    $\Gr$ into $A$ independently with probability $\frac{1}{2}$, then the probability of the event
  \begin{equation}
  \label{large_dev_k}
    |\sigma_A (X,{y_j})| \ge \eps\quad (j=1,\dots,k)
\end{equation}
is at most
$$\exp\left(\frac{-\eps^2 kn}2\right).$$
\label{l:M1}
\end{lemma}

{\bf Proof.} Denote by $\H_i$ ($i=0,\dots,k$) the event
$$
    |\sigma_A (X,{y_j})| \ge \eps\quad (j=1,\dots,i).
$$
We will prove by induction on $i$ that the probability
$\mathbb{P}_i$ of the event $\H_i$ is at most
$$\exp\left(\frac{-\eps^2 in}2\right).$$
The claim is obvious for $i=0$. 
%IS
%We prove that for any $i=1,\dots,k$ it is true for $i$ provided that it holds for $i-1$.
Now we prove that it is true for each $i=1,...,k$ whenever it holds for $i-1$.

Let
$$X_i' = \{x\in X:\,\exists j\in \{1,\dots,i-1\}:\, x+y_i\in X+y_j\},
\quad X_i= X\setminus X_i'.$$
By (\ref{small_inters}) we have
$$|X_i'|\le\eps n.$$
Therefore,
$$
\left|\sum_{x\in X_i'} A(x+y_i) - \frac{|X_i'|}{2}\right|
\le\frac{\eps n}2.
$$
Assuming that $\H_i$ holds, we see that
$$\left|\sum_{x\in X} A(x+y_i) - \frac{|X|}{2}\right|
\ge\eps n.
$$
Hence,
$$
\left|\sum_{x\in X_i} A(x+y_i) - \frac{|X_i|}{2}\right|
\ge \left|\sum_{x\in X} A(x+y_i) - \frac{|X|}{2}\right|
$$
$$
- \left|\sum_{x\in X_i'} A(x+y_i) - \frac{|X_i'|}{2}\right|
\ge \frac{\eps n}2.
$$
Thus, denoting by $\H_i'$ the event
$$
\left|\sum_{x\in X_i} A(x+y_i) - \frac{|X_i|}{2}\right| \ge \frac{\eps n}2,$$
we conclude that $\H_i'$ holds if $\H_i$ holds. Since for $x\in X_i$
the element $x+y_i$ does not belong to the sets $X+y_j,\,j<i,$
the event $\H_i'$ is independent of the events $\H_1,\dots,\H_{i-1}$.
Therefore, if $\mathbb{P}_i'$ is the probability of the event $\H_i'$,
then
\begin{equation}
\label{est_prob}
\mathbb{P}_i\le \mathbb{P}_{i-1}\mathbb{P}_i'.
\end{equation}

By Proposition 3 from \cite{M} (Hoeffding's theorem) we get
$$\mathbb{P}_i' \le \exp\left( -\frac12\left(\frac{\eps n/2}{\sqrt{|X_i|}/2}\right)^2\right)
\le \exp\left(-\frac{\eps^2n}2 \right).$$ Plugging in this estimate into
(\ref{est_prob}) and using the induction hypothesis we complete
the proof of the lemma. $\hfill\Box$

\begin{cor}
    Let $\mathbf{G}$ be a finite abelian group, and write $N = |\mathbf{G}|$.
    Fix $X\subset\Gr$ with $|X|=n$. Let $\eps \in (0,1/2]$ and $k\in\N$. Then
the probability that there exist $y_1,\dots,y_k\in \Gr$ satisfying
(\ref{small_inters}) and (\ref{large_dev_k}) is at most
$$\left(N\exp\left(\frac{-\eps^2 n}2\right)\right)^k.$$
\label{c:M2}
\end{cor}

\begin{cor}
    Let $\mathbf{G}$ be a finite abelian group, and write $N = |\mathbf{G}|$.
    Fix $X\subset\Gr$ with $|X|=n$. Let $\eps \in (0,1/2]$ and $k\in\N$. If
$$n\ge \frac{4\log N}{\eps^2},$$
then the probability that there exist $y_1,\dots,y_k\in \Gr$ satisfying
(\ref{small_inters}) and (\ref{large_dev_k}) is at most
$$\exp\left(\frac{-\eps^2 nk}{4}\right).$$
\label{c:M3}
\end{cor}

Now we will show that if $Y\subset\Gr$ is a large set and $\sigma_A (X,Y)$
is also large then for an appropriate $\eps$ there are many elements $y_1,\dots,y_k$
satisfying (\ref{small_inters}) and (\ref{large_dev_k}). Observe that if
we even do not assume that $y_1,\dots,y_k$ are distinct, this would follow
from (\ref{small_inters}).

%We denote
%$$X\circ X(z) = \sum_{x\in X} X(x+z).$$
%Equivalently, $X\circ X(z)$ is the number of representations of $z$
%as $x_1-x_2$ with $x_1,x_2\in X$. Clearly,
%\begin{equation}
%\label{sum_circ}
%\sum_{z\in\Gr}X\circ X(z) = |X|^2.
%\end{equation}
%Also, observe that for $y,y'\in\Gr$ we have
%\begin{equation}
%\label{inters_circ}
%|(X+y)\cap (X+y')| = X\circ X(y-y').
%\end{equation}

\begin{lemma}
	Let $\mathbf{G}$ be a finite abelian group,
	$X,Y\subset\Gr$ with $|X|=n$, and let $\eps \in (0,1/2]$. If
	$$\E(X,Y) = |X|^2|Y|/K,$$
	then for some
	$$k>\eps^2|Y|K/n$$
	there are $y_1,\dots,y_k\in Y$ satisfying condition
	(\ref{small_inters}).
	%\begin{equation}
	%\label{disc2}
	%(X+y_i)\cap\left(\bigcup_{j=1}^{i-1}(X+y_j)\right)\le\eps n/2
	%    \quad(i=2,\dots,k).
	%\end{equation}
	%and
	%$$
	%    |\sigma_A (X,{y_j})| \ge \eps/2\quad (j=1,\dots,k).
	%$$
	\label{l:disc2}
\end{lemma}

{\bf Proof.}
Let $\{y_1,\dots,y_k\}$ be the maximal subset of $Y$ satisfying (\ref{small_inters}). Denote
$$Z=\bigcup_{i=1}^k (X+{y_i}).$$
For any $z\in Z$ we denote by $f(z)$ the number of solutions of the equation
$$x+y=z,\quad x\in X,\,y\in Y.$$
By the choice of $k$, for any $y\in Y$ there are more than $\eps n$
values $x\in X$ such that $x+y\in Z$. Hence,
$$\sum_{z\in Z} f(z)> \eps n|Y|.$$
By the Cauchy--Schwartz inequality
$$\sum_{z\in Z}f^2(z) \ge |Z|^{-1}\left(\sum_{z\in Z} f(z)\right)^2
>\frac{\eps^2n|Y|^2}{k}.$$
Since
$$\E(X,Y) = \sum_{z\in\Gr} f^2(z),$$
we conclude that
$$n^2|Y|/K>\frac{\eps^2n|Y|^2}{k}$$
implying the required inequality for $k$.
%$$k>\eps^2|Y|K/(4n).$$

$\hfill\Box$

\begin{lemma}
    Let $\mathbf{G}$ be a finite abelian group,
    $X,Y\subset\Gr$, and let $\eps \in (0,1/2]$. 
%IS
%    Also, let $A$ be a random subset of $\mathbf{G}$ obtained by putting every element of $\Gr$ into $A$ independently with probability $\frac{1}{2}$. 
	Also, let $A \subset \mathbf{G}$ be a set.  
    If
$$|\sigma_A (X,Y)| \ge\eps,$$
then there is a set $Y'\subset Y$ such that $|Y'|\ge\eps|Y|$ and
any $y\in Y'$ satisfies the condition
$$|\sigma_A (X,{y})| \ge\eps/2.$$
\label{l:subset}
\end{lemma}

{\bf Proof.}
Denote
$$Y' = \{y\in Y:\,|\sigma_A (X,{y})|\ge\eps/2\},\quad Y''=Y\setminus Y'.$$
We have
$$
\left|\sum_{x\in X, y\in Y''} A(x+y) - \frac{|X|||Y''|}{2}\right|
\le\sum_{y\in Y''}\left|\sum_{x\in X} A(x+y) - \frac{|X|}{2}\right|
$$
$$
\le\sum_{y\in Y''}\eps|X|/2 \le \eps|X||Y|/2.
$$
Therefore,
$$
\left|\sum_{x\in X, y\in Y'} A(x+y) - \frac{|X|||Y'|}{2}\right|
\ge \left|\sum_{x\in X, y\in Y} A(x+y) - \frac{|X|||Y|}{2}\right|
$$
$$
- \left|\sum_{x\in X, y\in Y''} A(x+y) - \frac{|X|||Y''|}{2}\right|
\ge\eps|X||Y|/2.
$$
On the other hand,
$$
\left|\sum_{x\in X, y\in Y'} A(x+y) - \frac{|X|||Y'|}{2}\right|
\le|X||Y'|/2.$$
Thus, $|Y'|\ge\eps|Y|$ as required.
$\hfill\Box$

\bigskip

Combining Lemmas 
\ref{l:subset}  
and 
%\ref{l:discr}, 
\ref{l:disc2}, we get the following
corollary.

\begin{cor}
    Let $\mathbf{G}$ be a finite abelian group, %, and write $N = |\mathbf{G}|$.
    $X,Y\subset\Gr$ with $|X|=n$, and let $\eps \in (0,1/2]$. 
%IS
%    Also, let $A$ be a random subset of $\mathbf{G}$ obtained by putting every element of $\Gr$ into $A$ independently with probability $\frac{1}{2}$. 
Also, let $A \subset \mathbf{G}$ be a set. 
    If
$$\E(X,Y) \le |X|^2|Y|/K,\quad |\sigma_A (X,Y)| \ge\eps,$$
then for some
%$$k> \eps^2|Y|/(2n)$$
$$k>\eps^4|Y|K/(4n)$$
there are $y_1,\dots,y_k\in Y$ satisfying conditions
\begin{equation}
\label{disc2}
\left|(X+y_i)\cap\left(\bigcup_{j=1}^{i-1}(X+y_j)\right)\right|\le\eps n/2
    \quad(i=2,\dots,k) 
\end{equation}
%$$
% \sum_{j=1}^{i-1}\left(|X+{y_i}\cap X+{y_j}|\right) \le \eps n/2
%    \quad(i=2,\dots,k)
%$$
and the condition
$$
    |\sigma_A (X,{y_j})| \ge \eps/2\quad (j=1,\dots,k).
$$
\label{c:disc}
\end{cor}
{\bf Proof.}
We take a subset $Y'\subset Y$, in accordance
with Lemma \ref{l:subset}. Let
$$\E(X,Y') = |X|^2|Y'|/K'.$$
Since $|Y'|\ge\eps|Y|$, $\E(X,Y')\le \E(X,Y)\le|X|^2|Y|/K$, we have
$$K'\ge\eps K.$$
Applying Lemma \ref{l:disc2} (with $\eps/2$ instead of $\eps$) to the set $Y'$ we get desired $y_1,\dots,y_k\in Y'$
with
$$k>\eps^2|Y'|K'/(4n)\ge\eps^4|Y|K/(4n).$$
$\hfill\Box$

%Corollaries \ref{c:disc} and \ref{c:M3} immediately imply the main result of this section.

%\begin{proposition}
%    Let $\mathbf{G}$ be a finite abelian group, and write $N = |\mathbf{G}|$.
%    Fix $X\subset\Gr$ with $|X|=n$, $\eps \in (0,1/2]$, and $m\in\N$. Also, let
%$A$ be a random subset of $\mathbf{G}$ obtained by putting every element of
%    $\Gr$ into $A$ independently with probability $\frac{1}{2}$. If
%$$n\ge \frac{4\log N}{\eps^2},$$
%then the probability that there exist $Y\subset \Gr$ satisfying
%$$|Y|\ge m,\quad |\sigma_A (X,Y)| \ge\eps,$$
%is at most
%$$\exp\left(\frac{-\eps^4 m}{32}\right).$$
%\label{p:est_prob}
%\end{proposition}

Corollaries  
%\ref{c:disc2} 
\ref{c:M3} 
and 
\ref{c:disc}
immediately imply
the main result of this section.

\begin{proposition}
	Let $\mathbf{G}$ be a finite abelian group, and write $N = |\mathbf{G}|$.
	Fix $X\subset\Gr$ with $|X|=n$, $\eps \in (0,1/2]$, $K\ge1$, and $m\in\N$. 
	Let $A$ be a random subset of $\mathbf{G}$ obtained by putting every element of
	$\Gr$ into $A$ independently with probability $\frac{1}{2}$. If
	$$n\ge \frac{4\log N}{\eps^2},$$
	then the probability that there exist $Y\subset \Gr$ satisfying
	$$|Y|\ge m,\quad |\sigma_A (X,Y)| \ge\eps,\quad \E(X,Y)\le|X|^2|Y|/K,$$
	is at most
	$$\exp\left(\frac{-\eps^6 mK}{64}\right).$$
	\label{p:est_prob2}
\end{proposition}

\bigskip
\section{The proof of the main result}
\label{sec:improvement}
\bigskip

%{\bf Proof of Theorem \ref{t:main3}.}

In this section we obtain our main Theorem \ref{t:main3}.\\
Let $N=|\Gr|$. Without loss of generality we assume that
\begin{equation}
\label{w_suppos}
w(N)\le\log\log(N+3).
\end{equation}
Denote
\begin{equation}
\label{choice_eps}
w_1(N) =\sqrt{w(N)},\quad \eps=w(N)^{-1/13}
\end{equation}
and
%IS n_0 =... -> n_0:= ...
%SK remove ":" similarly to definitions of W_1, \eps, \tilde n_1
$$\tilde n_0 =\left[w_1(N)\log N (\log \log N)^2\right],
\quad \tilde n_1=2\tilde n_0.$$
Assume that
    $$|\tilde X|\ge w(N)\log N (\log \log N)^2,\quad
    |Y| \ge m:=w(N) \log N(\log \log N)^{10} \,,
    $$
$$|Y|\ge|\tilde X|,\quad|\sigma_A (\tilde X,Y)|\ge\eps/2,$$
and $w(N)$ is large enough as well as $N$. We have to prove that
the probability of existence of such sets $\tilde X,Y$ is small.
In view of Theorem \ref{t:Mrazovic}, we can assume that $|\tilde X|, |Y| \ \ll \log^{5/2} N$, say, 
%IS
because 
otherwise the probability of  the event $|\sigma_A (\tilde X,Y)| \ge \eps/2$ is $o(1)$.
Since $|\tilde X|\ge \tilde n_1$, we can split $\tilde X$ into sets $X_i$ with
$\tilde n_0\le |X_i|\le \tilde n_1$ in an arbitrary way.
For some $i$ we have $|\sigma_A (X_i,Y)| \ge |\sigma_A (\tilde X,Y)|\ge\eps/2$.
Take $X = X_i$. If these sets $X,Y$ exist, then the following event $\H_0$ happens:\newline
there exist $X,Y\subset\Gr$ such that
$$|X|\le 2w_1(N) \log N(\log\log N)^2,\quad |Y|\ge m,$$
%IS
%According to Misha's (26) remark, I have replaced all \eps -> \tilde{\eps}=\eps/4 below.
$$
	\left|\sum_{x\in X} \sum_{y\in Y} \left(A(x+y) - \frac{1}{2}\right)\right| 
	\ge 
	\tilde{\eps} w_1(N) (\log N)(\log\log N)^2|Y| \,,
$$
where $\tilde{\eps} = \eps/4$. 
Our aim is to prove that the probability $\P_0$ of the event
$\H_0$ tends to $0$ as $N\to\infty$. We will consider the
family of events $\{\H_j\}$, $j\ge0$. We say that $\H_j$ happens if
there exist $X,Y\subset\Gr$ such that
\begin{equation}
\label{sizeX} |X|\le 2w_1(N) \log N(\log\log N)^2,
\end{equation}
%IS 
\begin{equation}
%\label{sizeY} |Y|\ge m:= w(N) \log N(\log\log N)^{11},
\label{sizeY} |Y|\ge m= w(N) \log N(\log\log N)^{10},
\end{equation}
\begin{equation}
\label{sizeE}
\E(X,Y) \le |X|^2|Y|10^{-j},
\end{equation}
\begin{equation*}\label{sizesigma} 
\left|\sum_{x\in X} \sum_{y\in Y} \left(A(x+y) -
\frac{1}{2}\right)\right| 
\ge 
\end{equation*}
\begin{equation}\label{sizesigma}
\ge 
\left(1-\frac j{\log\log N}\right) \tilde{\eps} w_1(N) (\log N)(\log\log N)^2|Y|,
\end{equation}
where $j\ge0$. We denote by $\P_j$ the probability of the event $\H_j$.
Let $j_0=[(\log\log N)/2]$. We observe that, due to (\ref{w_suppos}) 
and (\ref{sizeX}),
$$10^{j_0}>|X|.$$ 
%IS
%SK
Hence, $\P_{j_0}=0$ 
(because (\ref{sizeE}) does not hold  for $j=j_0$ due to $\E(X,Y) \ge |X||Y|$), 
and we will consider that $j<j_0$.

We will estimate $\P_j$ in terms of $\P_{j+1}$. Let $X,Y\subset\Gr$
%IS
%satisfying 
satisfy 
(\ref{sizeX})--(\ref{sizesigma}). Denote
$$K = \frac{|X|^2 |Y|}{\E(X,Y)},\quad M = 10^{j+1}.$$
If $M\le K$ then the event $\H_{j+1}$ holds. Thus we can assume
that $M\ge K$. Applying Corollary \ref{c:SS_dim} to the sets $X$,
$Y$, we find $X',X''\subset X$ such that $X=X' \bigsqcup X''$,
$\dim(X') \ll M (\log \log N)^2\ll 10^{j+1}(\log\log N)^2$,
$\E(X',Y) \gg \E(X,Y)$ and $\E(X'',Y)\le |Y| |X''|^2/M$.

%IS
%Consider first the case where inequality
Firstly, consider the case where  the inequality  
\begin{equation}\label{f:supposX'2}
\left|\sum_{x\in X'} \sum_{y\in Y} \left(A(x+y) -
\frac{1}{2}\right)\right| \ge \tilde{\eps} w_1(N) (\log N)(\log\log N)|Y|
\end{equation}
does not hold. Then
$$
\left|\sum_{x\in X''} \sum_{y\in Y} \left(A(x+y) - \frac{1}{2}\right)\right|
	\ge 
$$
$$
\ge\left|\sum_{x\in X} \sum_{y\in Y} \left(A(x+y) - \frac{1}{2}\right)\right|
%$$
%$$
- \left|\sum_{x\in X'} \sum_{y\in Y} \left(A(x+y) -
\frac{1}{2}\right)\right| 
\ge
$$
$$
\ge \left(1-\frac {j+1}{\log\log N}\right) \tilde{\eps} w_1(N) (\log N)(\log\log N)^2|Y|,
$$
and (\ref{sizeX})--(\ref{sizesigma}) hold for $j+1$ instead of $j$
and $X''$ instead of $X$. Again, $\H_{j+1}$ holds.

Now consider the case where (\ref{f:supposX'2}) holds.
Then we have
%IS n_0:= ... -> n_0= ....
$$|X'| \ge \tilde{\eps} w_1(N) (\log N)(\log\log N) := n_0 .$$
%=\left[\eps w_1(N) (\log N)(\log\log N)\right].$$
Let
$$n_\nu =2^{\nu} n_0,\quad \nu\le\nu_0,$$
where $\nu_0$ is defined by 
%IS The sign in the ineq. below is changed
$$n_{\nu_0} \le 2w_1(N) \log N(\log\log N)^2 < n_{\nu_0+1}.$$
%IS 
Clearly, $\nu_0 \ll \log \log \log N+ \log (1/\eps)$. 
By Lemma \ref{l:count_sm_dim}, the number of such sets $X'$
with $|X'|=n$, $n_\nu\le n< n_{\nu+1}$, is at most
\begin{equation}\label{tmp:13.09_1}
e^{Cn_\nu 10^j(\log\log N)^2}\le e^{2C 10^jw_1(N)\log N(\log\log N)^4} \,,
\end{equation}
where $C$ is an absolute constant.

Next, for these sets $X'$
inequality (\ref{f:supposX'2}) implies that
%IS = -> \ge
$$|\sigma_A(X',Y)| \ge \eps'= \tilde{\eps} w_1(N) (\log N)(\log\log N)/n_{\nu+1}
%=
\ge
$$
$$ 
\ge
\tilde{\eps} w_1(N) (\log N)(\log\log N)/(2n_\nu).$$ 
We have
$$n_\nu(\eps')^2/(4\log N)\ge n_\nu \tilde{\eps}^2 w_1(N)^2(\log N)^2(\log\log N)^2/(16n_\nu^2\log N)$$
$$=\tilde{\eps}^2 w_1(N)^2(\log N)(\log\log N)^2/(16n_\nu)\ge \tilde{\eps}^2 w_1(N)/32 >1,$$
where we have used (\ref{choice_eps}), and we are in position to use
Proposition \ref{p:est_prob2}. We have
$$\E(X',Y) \le \E(X,Y) \le |X|^2|Y|/10^j $$
$$\le (2w_1(N)\log N(\log\log N)^2)^2 |Y|/10^j = n_\nu^2|Y|/K' \le |X'|^2|Y|/K',$$
where
$$K' = 10^j \left(n_\nu/(2w_1(N)\log N(\log\log N)^2)\right)^2$$ 
$$= 10^j n_\nu^2/\left(4w_1(N)^2(\log N)^2(\log\log N)^4\right).$$
Thus, for any such $X'$ the probability $\P_{j,\nu}(X')$
of the existence of a set $Y$ satisfying this inequality is at most
$$\exp\left(\frac{-(\eps')^6 m \max\{K',1\}}{64}\right) \le \exp\left(\frac{-(\eps')^6 m K'\}}{256}\right).$$
%Hence,
We have 
$$
(\eps')^6 K' = 10^j\tilde{\eps}^6 w_1(N)^4(n_\nu)^{-4}(\log N)^4(\log\log N)^2/256
\ge 
$$
\begin{equation}\label{tmp:13.09_2}
\ge 10^j\eps^6 (\log\log N)^{-6}/ 2^{18}.
\end{equation}
Next, the probability $\P_{j,\nu}$ of the existence of a set $X'$,
$|X'|=n$, $n_\nu\le |X'|< n_{\nu+1}$, that can be obtained by our
%IS References added. 
construction, is bounded by (see (\ref{tmp:13.09_1}), (\ref{tmp:13.09_2}))
$$\exp\left(2Cw_1(N)10^j\log N(\log\log N)^4 - 10^j\eps^6(\log\log N)^{-6}m/2^{26} \right).$$
Taking into account (\ref{choice_eps}) 
%IS 
and the definition of  the parameter $m$, we get
$$\P_{j,\nu} \le \exp\left(-10^jw_1(N)\log N(\log\log N)^4\right).$$
Taking the sum over $\nu$, we find
%IS \mu -> \nu
$$\P_j\le \P_{j+1} +\sum_{\nu} \P_{j,\nu} \le \P_{j+1} +
\exp\left(-10^jw_1(N)\log N(\log\log N)^4/2\right).$$
Finally,
$$\P_0 \le\sum_{j=0}^{j_0-1}\exp\left(-10^jw_1(N)\log N(\log\log N)^4/2\right)$$
$$\le\exp\left(-w_1(N)\log N(\log\log N)^4/3\right).$$
$\hfill\Box$

\bigskip
\section{Appendix}
\label{sec:appendix}
\bigskip

In this section we prove Proposition \ref{p:M}.

\bigskip

By sections 3, 4 of \cite{M} there are sets $S \subset X$, $T\subset Y$, $s=|S|$, $t=|T|$
such that
\begin{equation}\label{f:app1}
\E(S,T) \le 2st + \frac{2s^2 t^2}{|X|^2 |Y|^2} \cdot \E(X,Y) \,,
\end{equation}
and
\begin{equation}\label{f:app2}
|\sigma_A (X,Y) - \sigma_A (S,T)| \le 6\sqrt\frac{|Y|}{st} \,.
\end{equation}
Here $s,t$ are parameters and we choose $s=\frac{2000 \log N}{\eps^4} \le |X|$ and $t=\frac{K|Y|\eps^2}{10 \log N}$.
The left--hand side of (\ref{f:app2}) is less than $\eps/2$.
On the other hand, by the large deviations low (see Proposition 3 and calculations after this proposition from \cite{M}) and (\ref{f:app1}) the probability $\mathbb{P}$ of  $|\sigma_A (S,T)| \ge \eps/2$ is bounded by
$$
\mathbb{P} \ll \exp\left(-\frac{\eps^2 s^2 t^2}{2\E(S,T)} \right)
\ll
\exp\left( - \min\left\{ \frac{\eps^2 st}{8}, \frac{\eps^2 |X|^2 |Y|^2}{8 \E(X,Y)} \right\} \right)
\ll
$$
$$
\ll
\exp \left( - \min\left\{ 25 K|Y|, \frac{\eps^2 K|Y|}{8} \right\} \right)
\ll
\exp \left( - \frac{\eps^2 K|Y|}{8} \right) \,.
$$
Thus, the final probability (\ref{f:p_M})
does not exceed
$$
N^{s+t} \exp \left( - \frac{\eps^2 K|Y|}{8} \right)
\ll
\exp\left( \frac{2000 \log^2 N}{\eps^4} + \frac{K|Y|\eps^2}{10}  - \frac{\eps^2 K|Y|}{8} \right)
\le
$$
$$
\le
\exp\left( \frac{2000 \log^2 N}{\eps^4} - \frac{\eps^2 rK}{40} \right) %\,.
$$
as required.

{}
\bigskip
%\newpage

\noindent{S.V.~Konyagin\\
Steklov Mathematical Institute of Russian Academy of Sciences, Moscow

{\tt konyagin@mi.ras.ru}

\bigskip

\noindent{I.D.~Shkredov\\
Steklov Mathematical Institute of Russian Academy of Sciences, Moscow

{\tt ilya.shkredov@gmail.com}

\end{document}